\documentclass{amsart}

\usepackage{amsmath}
\usepackage{amssymb}
\usepackage{amsthm}
\usepackage{tikz}
\usepackage{graphicx}
\usetikzlibrary{decorations.pathreplacing}

\usepackage{hyperref}

\newcommand{\rsp}{\raisebox{0em}[2.7ex][1.3ex]{\rule{0em}{2ex} }}

\newtheorem{lem}{Lemma}
\newtheorem{prop}[lem]{Proposition}
\newtheorem{thm}[lem]{Theorem}

\title{On the Babylonian Division of Trapezoids}
\author{F. Lemmermeyer}

\begin{document}

\maketitle

One of the most remarkable achievements\footnote{For a detailed
  description cf. \cite{BBS}; see also \cite{LemMI}, \cite{DamPyth}
  and \cite{Lem4000}.} of Babylonian mathematics is the bisection of
trapezoids, which can be traced back to the time of Sargon of Akkad.
On tablet\footnote{This tablet is discussed e.g. in \cite{Frib,FribT}
  and \cite{GoncTH}. IM stands for the Iraq Museum in Baghdad; how
  many of these tablets survived the looting during the Iraq War is
  probably still unknown today; see \cite{Polk}.}  IM~58045 one finds
the drawing of a trapezoid of height $h=12$, whose parallel sides have
lengths $a = 17$ and $c = 7$. The problem apparently consisted in
dividing the trapezoid by a line parallel to $a$ and $c$ into two
trapezoids of equal area and computing their side lengths. Conditions
such as the parallelism of the transversal are rarely stated explicitly,
but follow from the given solutions. For finding solutions
in integers or rational numbers, one must solve a Diophantine
equation. In fact, the segment of length $b$ that is parallel to the
parallel sides of lengths $a$ and $c$ of a trapezoid divides the
trapezoid into two partial trapezoids of equal area if the relation
$a^2 + c^2 = 2b^2$ holds.

\begin{thm}
  Given a trapezoid with parallel sides $a$ and $c$, the transversal
  parallel to these sides bisects the trapezoid into two parts of equal
  area if it has length $b$, where
  \begin{equation}\label{BabT} a^2 + c^2 = 2b^2. \end{equation}
\end{thm}

Guessing a relation such as (\ref{BabT}) without some form of
justification is virtually impossible. Nobody has ever discovered
the Pythagorean theorem without a proof (or, at least, a heuristic
explanation). In this article\footnote{This is based on \cite{Lem4000},
  where I screwed up the description of the ``Babylonian Pythagoras''.}
I will present a few proofs that are based on diagrams believed to
be familiar to the Babylonian scribes.

\begin{center}
  \begin{tikzpicture}[scale=0.4]
    \fill[green, opacity=0.3] (0,0) -- (17,0) -- (15,2) -- (2,2) -- cycle;
    \fill[red,   opacity=0.3] (2,2) -- (15,2) -- (12,5) -- (5,5) -- cycle;
    \draw (0,0) -- (17,0) -- (12,5) -- (5,5) -- cycle;
    \draw (2,2) -- (15,2);
  \draw [decorate, decoration = {brace}]
        (17,-0.2) --  (0,-0.2);
  \draw (8.5,-0.3) node[anchor=north] {$a$};
  \draw [decorate, decoration = {brace}]
        (14.95,1.8) --  (2.05,1.8);
  \draw (8.5,1.7) node[anchor=north] {$b$};
  \draw [decorate, decoration = {brace}]
        (5,5.2) --  (12,5.2);
  \draw (8.5,5.3) node[anchor=south] {$c$};
  \end{tikzpicture}  
\end{center}

The trapezoid with sides $a = 17$, $b = 13$ and height $2$
has the same area $30$ as the trapezoid with sides $b = 13$, $c = 7$
and height $3$, and both have the same interior angle $45^\circ$.
There are various methods for deriving the Babylonian Theorem
$a^2 + c^2 = 2b^2$, and we will present three of them.

\section{Euclid}\label{S1}

Euclid is best known for his {\em Elements}; his book on the 
Division of Figures, which is preserved only partially in
an Arabic recension, is less well known. In Prop.~4 of this work,
the halving of trapezoids is discussed.

The essential idea behind Euclid's derivation\footnote{See
  \cite[p.~35]{ED}.} of the relation $a^2 + c^2 = 2b^2$ for the
bisection of trapezoids is the observation that the areas of similar
triangles with bases $x$ and $y$ are in the ratio $x^2 : y^2$. This is
a special case of Proposition 23 in the sixth book of Euclid's {\em Elements}.

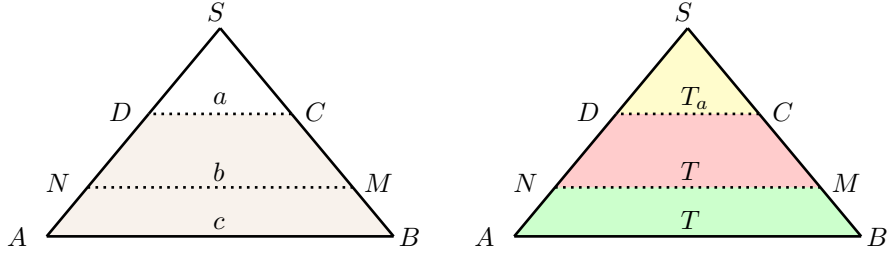
\begin{figure}[ht!]
\begin{center}
\begin{tikzpicture}[scale=0.27]
\fill[brown, opacity=0.1]
  (5,6) -- (0,0) -- (17,0) -- (12,6) -- cycle;
\draw [line width=1pt] (5,6)-- (0,0);
\draw [line width=1pt] (0,0)-- (17,0);
\draw [line width=1pt] (17,0)-- (12,6);
\draw [line width=1pt, dotted] (12,6)-- (5,6);
\draw [line width=1pt, dotted] (2,2.4)-- (15,2.4);
\draw [line width=1pt] (8.5,10.2)-- (12,6);
\draw [line width=1pt] (8.5,10.2)-- (5,6);
\draw (-2.4,0.9) node[anchor=north west] {$A$};
\draw (16.8,0.9) node[anchor=north west] {$B$};
\draw (-0.5,3.5) node[anchor=north west] {$N$};
\draw (15.1,3.5) node[anchor=north west] {$M$};
\draw (2.6,7) node[anchor=north west] {$D$};
\draw (12.2,7) node[anchor=north west] {$C$};
\draw (7.4,11.9) node[anchor=north west] {$S$};
\draw (7.7,1.4) node[anchor=north west] {$c$};
\draw (7.7,4.1) node[anchor=north west] {$b$};
\draw (7.7,7.5) node[anchor=north west] {$a$};
\end{tikzpicture} \quad \begin{tikzpicture}[scale=0.27]
\fill[green,  opacity=0.2] (0,0) -- (17,0) -- (15,2.4) -- (2,2.4) -- cycle;
\fill[red,    opacity=0.2] (12,6) -- (5,6) -- (2,2.4) -- (15,2.4) -- cycle;
\fill[yellow, opacity=0.2] (5,6) -- (12,6) -- (8.5,10.2) -- cycle;
\draw [line width=1pt] (5,6)-- (0,0);
\draw [line width=1pt] (0,0)-- (17,0);
\draw [line width=1pt] (17,0)-- (12,6);
\draw [line width=1pt, dotted] (12,6)-- (5,6);
\draw [line width=1pt, dotted] (2,2.4)-- (15,2.4);
\draw [line width=1pt] (8.5,10.2)-- (12,6);
\draw [line width=1pt] (8.5,10.2)-- (5,6);
\draw (-2.4,0.9) node[anchor=north west] {$A$};
\draw (16.8,0.9) node[anchor=north west] {$B$};
\draw (-0.5,3.5) node[anchor=north west] {$N$};
\draw (15.1,3.5) node[anchor=north west] {$M$};
\draw (2.6,7) node[anchor=north west] {$D$};
\draw (12.2,7) node[anchor=north west] {$C$};
\draw (7.4,11.9) node[anchor=north west] {$S$};
\draw (7.7,1.7) node[anchor=north west] {$T$};
\draw (7.7,4.1) node[anchor=north west] {$T$};
\draw (7.7,7.8) node[anchor=north west] {$T_a$};
\end{tikzpicture}
\end{center}
\caption{Division of trapezoids in Euclid}\label{AbbTTE}
\end{figure}

Assume that $ABCD$ is a trapezoid with parallel sides $AB$ and $CD$, and
that it is bisected by the transversal $MN$ parallel to $AB$
(see Fig.~\ref{AbbTTE}). We complete the trapezoid to triangle $ABS$ as
indicated. Denote the areas of the partial triangles $SCD$, $SMN$, and
$SAB$ by $T_a$, $T_b$, and $T_c$, respectively, and let $T$ be the area
of the halved trapezoid. Then the following relations hold:
$$ T_a + T = T_b, \quad T_c - T = T_b, \quad \text{hence} \quad
    T_a+T_c = 2T_b.$$

Since the areas $T_a$, $T_b$, and $T_c$ are to one another as
$a^2 : b^2 : c^2$, it follows, almost without any calculation,
that $a^2 + c^2 = 2b^2$.

\begin{center}
  \begin{tikzpicture}[scale=0.6]  
  \fill[green, opacity=0.1] (0,0) -- (0,7) -- (7,7) --
    (7,0) -- (6,1) -- (6,6) -- (1,6) -- (1,1) -- cycle;
  \fill[green, opacity=0.3] (0,0) --  (7,0) -- (6,1) -- (1,1) -- cycle;
  \fill[red,   opacity=0.1] (1,1) --  (1,6) -- (6,6) -- (6,1) --
     (4,3) -- (4,4) -- (3,4) -- (3,3) -- cycle;
  \fill[red,   opacity=0.3] (1,1) --  (6,1) -- (4,3) -- (3,3) -- cycle;    
  \draw (0,0) -- (7,0) -- (7,7) -- (0,7) -- cycle;
  \draw (1,1) -- (6,1) -- (6,6) -- (1,6) -- cycle;
  \draw (3,3) -- (4,3) -- (4,4) -- (3,4) -- cycle;
  \draw (0,0) -- (7,7);
  \draw (7,0) -- (0,7);
  \draw [decorate, decoration = {brace}]
        (7,-0.2) --  (0,-0.2);
  \draw (3.5,-0.3) node[anchor=north] {$a$};
  \draw [decorate, decoration = {brace}]
        (5.95,0.8) --  (1.05,0.8);
  \draw (3.5,0.7) node[anchor=north] {$b$};
  \draw [decorate, decoration = {brace}]
        (3.95,2.9) --  (3.05,2.9);
  \draw (3.5,2.7) node[anchor=north] {$c$};
  \end{tikzpicture}
\end{center}

This proof becomes completely trivial in the special cases where
the base angle of the trapezoid is $45^\circ$ (to which the general
case may be reduced by using similarity of triangles). In this case,
four trapezoids form a square, and a look at the figure immediately
tells us that if $b$ bisects the trapezoid with side lengths $a$ and $c$,
then $a^2 - b^2 = b^2 - c^2$, which in turn is equivalent to
$a^2 + c^2 = 2b^2$.

\section{Babylonian and Pythagorean Triples}

Below, we shall call a triple of natural numbers $(a, b, c)$
with $a^2 + c^2 = 2b^2$ a Babylonian triple.
If $(x, y, z)$ is a Pythagorean triple, i.e., a solution of the equation
$x^2 + y^2 = z^2$ in natural numbers, and if $x < y$, then
$(y+x)^2 + (y-x)^2 = 2z^2$ and therefore $(a, b, c) = (y+x, z, y-x)$
is a Babylonian triple.

From the Pythagorean rule for the construction of right-angled
triangles in numbers, i.e., the triples
$$ (x, y, z) = \Big(2m+1, \frac{(2m+1)^2-1}2, \frac{(2m+1)^2+1}2\Big), $$
we thus obtain one of the two families of Babylonian triples $(a, b, c)$
in Tab.~\ref{TabT}.

If, on the other hand, one uses Plato's family of Pythagorean triples, i.e.,
$$ (x, y, z) = (2m, m^2-1, m^2+1), $$
then correspondingly, for all even values of $m$, one obtains
a different family of Babylonian triples.

\begin{table}[ht!]
$$ \begin{array}{rrr|rrr}
  \rsp  x &  y &  z &  a &  b &  c \\ \hline
  \rsp  3 &  4 &  5 &  7 &  5 &  1 \\
  \rsp  5 & 12 & 13 & 17 & 13 &  7 \\
  \rsp  7 & 24 & 25 & 31 & 25 & 17 \\
  \rsp  9 & 40 & 41 & 49 & 41 & 31 \\
  \rsp 11 & 60 & 61 & 71 & 61 & 49  
  \end{array} \qquad \qquad
  \begin{array}{rrr|rrr}
  \rsp  x &  y &   z &   a &   b &   c \\ \hline
  \rsp  4 &  3 &   5 &   7 &   5 &   1 \\
  \rsp  8 & 15 &  17 &  23 &  17 &   7 \\
  \rsp 12 & 35 &  37 &  47 &  37 &  23 \\
  \rsp 16 & 63 &  65 &  79 &  65 &  47 \\
  \rsp 20 & 99 & 101 & 119 & 101 &  79 
\end{array} $$

\caption{Sequences of trapezoids of the first (left) and the second kind 
  (right)}\label{TabT}
\end{table}

It turns out that the Babylonians knew these trapezoids, as is
witnessed by a problem on tablet AO 17264, which was already discussed
by Neugebauer and can also be found in Gandz\footnote{See
  \cite[p.~126ff.]{Neugebauer} and \cite{Gandz}, as well as
  \cite{Caveing6}.  This tablet appears to date from the Kassite
  period between the 16th and 12th centuries BCE.}. A trapezoidal
field whose parallel sides have lengths $a = 51$ and $c = 213$, and
whose other side lengths are $135$ and $81$, is to be divided
by parallel strips among six brothers in such a way that the first two
receive the same share, likewise the third and fourth, as well as the
fifth and sixth.

Before we discuss possible derivations of the Babylonian identity
(\ref{BabT}) and the connection with Pythagorean triples, we
supply some background, in particular the Four Bricks diagram,
which is the basis of many procedures in Babylonian metric algebra.

\section{Four Bricks}

The most important figure  in Babylonian Algebra consists of four
bricks\footnote{See H\o{}yrup \cite{Hoy}, \cite{LemmBA}, and the many
  publications of Aldo Bonet \cite{Bonet}, as well as \cite{Cid} for a
  presentation of the Four Bricks diagram without its historical background.}. 
It may be used for proving the Theorem of Pythagoras, for solving quadratic 
equations, or for halving trapezoids.

\subsection{Some Identities}

\begin{figure}[ht!]
\begin{center}
  \begin{tikzpicture}[scale=0.4]
  \fill[green, opacity=0.3] (0,0) rectangle (5,5);
  \fill[red,   opacity=0.3] (0,5) rectangle (2,7);
  \fill[brown, opacity=0.3] (2,5) rectangle (7,7);
  \fill[brown, opacity=0.3] (5,0) rectangle (7,5);
  \draw[dotted] (0,5) -- (2,5);
  \draw (0,0) rectangle (7,7);
  \draw (5,0) -- (5,2) -- (0,2);
  \draw (5,2) -- (5,5) -- (7,5);
  \draw (5,5) -- (2,5) -- (2,7);
  \draw (2,5) -- (2,2);
  \draw [decorate, decoration = {brace}]  
      (0,7.1) --  (1.95,7.1);
  \draw (1,7.2) node[anchor=south] {$b$};
  \draw [decorate, decoration = {brace}]  
      (4.95,-0.1) --  (0,-0.1);
  \draw (2.5,-0.2) node[anchor=north] {$a$};
  \draw [decorate, decoration = {brace}]  
      (7,-0.1) --  (5.05,-0.1);
  \draw (6,-0.2) node[anchor=north] {$b$};
\end{tikzpicture} \qquad  
  \begin{tikzpicture}[scale=0.4]
  \fill [brown, opacity=0.3] (0,0) rectangle (7,7);
  \fill [white] (2,2) rectangle (5,5);
  \draw (0,0) rectangle (7,7);
  \draw (5,0) -- (5,2) -- (0,2);
  \draw (5,2) -- (5,5) -- (7,5);
  \draw (5,5) -- (2,5) -- (2,7);
  \draw (2,5) -- (2,2);
  \draw [decorate, decoration = {brace}]  
      (2.05,2.1) --  (4.95,2.1);
  \draw (3.5,2.2) node[anchor=south] {$a-b$};
  \draw [decorate, decoration = {brace}]  
      (0,7.1) --  (1.95,7.1);
  \draw (1,7.2) node[anchor=south] {$b$};
  \draw [decorate, decoration = {brace}]  
      (4.95,-0.1) --  (0,-0.1);
  \draw (2.5,-0.2) node[anchor=north] {$a$};
  \draw [decorate, decoration = {brace}]  
      (7,-0.1) --  (5.05,-0.1);
  \draw (6,-0.2) node[anchor=north] {$b$};
\end{tikzpicture} \qquad   \begin{tikzpicture}[scale=0.4]
  \fill[green, opacity=0.2] (0,0) rectangle (5,5);
  \fill[green, opacity=0.2] (2,2) rectangle (7,7);
  \fill[red,  opacity=0.2] (5,0) rectangle (7,2);
  \fill[red,  opacity=0.2] (0,5) rectangle (2,7);    
  \draw (0,0) rectangle (7,7);
  \draw (5,0) -- (5,5) -- (0,5);
  \draw (2,7) -- (2,2) -- (7,2);
  \draw [decorate, decoration = {brace}]  
      (2.05,2.1) --  (4.95,2.1);
  \draw (3.5,2.2) node[anchor=south] {$a-b$};
  \draw [decorate, decoration = {brace}]  
      (0,7.1) --  (1.95,7.1);
  \draw (1,7.2) node[anchor=south] {$b$};
  \draw [decorate, decoration = {brace}]  
      (4.95,-0.1) --  (0,-0.1);
  \draw (2.5,-0.2) node[anchor=north] {$a$};
  \draw [decorate, decoration = {brace}]  
      (7,-0.1) --  (5.05,-.1);
  \draw (6,-0.2) node[anchor=north] {$b$};
  \end{tikzpicture}
  \caption{Four Bricks}\label{Abb4B1}
\end{center}
\end{figure}
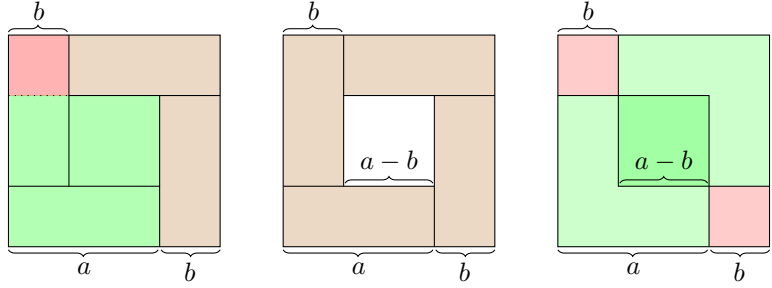

The Four Bricks diagram in the center of Fig.~\ref{Abb4B1} is the
source of various algebraic identities. It clearly shows that
$(a+b)^2 = (a-b)^2 + 4ab$, which is (\ref{E4B2}) below.

The left diagram in Figure~\ref{Abb4B1} similarly visualizes
$(a+b)^2 = a^2 + b^2 + 2ab$, which is (\ref{E4B1}) below; identity
(\ref{E4B1}) also follows from this diagram.
\begin{align}
  \label{E4B1}  (a+b)^2 & = a^2 + 2ab + b^2; \\
  \label{E4B2}  (a+b)^2 & = (a-b)^2 + 4ab; \\  
  \label{E4B3}  (a-b)^2 & = a^2 - 2ab + b^2.
\end{align}
We remark in passing that the geometric interpretation of the equation
$(a-b)^2 + 2ab = a^2 + b^2$  is Euclid's Proposition II.7, and that of
$(a+b)^2 + (a-b)^2 = 2(a^2 + b^2)$ is Prop.~II.10.

In the right diagram , $(a+b)^2$ is covered by the two green squares
with side lengths $a$ and the two red squares with side lengths $b$,
with the central deep green square with side length $a-b$ covered
twice. This implies
\begin{equation}\label{E4B4}
  (a+b)^2 + (a-b)^2 = 2a^2 + 2b^2.
\end{equation}

\subsection{The AGM Inequality}

Clearly, the area of the large square is at least as large as the
total area of the four rectangles, with equality if and only if the
central square vanishes, i.e., when $a = b$. This gives the inequality
$$ (a + b)^2 \ge 4ab, $$
and after dividing both sides by $4$ and taking the (positive)
square root we obtain the following

\begin{prop}\label{SAGM}
  For positive real numbers $a$ and $b$ we have
  $$ \frac{a + b}{2} \ge \sqrt{ab} $$
  with equality if and only if $a = b$.
\end{prop}

The identity (\ref{E4B4}) (as well as its geometric interpretation)
immediately implies $(a+b)^2 \le 2a^2 + 2b^2$, which in turn yields
the inequality between the quadratic and the arithmetic mean:

\begin{prop}
  For positive real numbers $a$ and $b$ we have
  $$ \frac{a+b}2 \le \sqrt{\frac{a^2 + b^2}2}. $$
\end{prop}

\subsection{Theorem of Pythagoras}

\begin{figure}[ht!]
\begin{center}
\begin{tikzpicture}[scale=0.7]
  \fill [brown, opacity = 0.3] (0,0) rectangle (4,2);
  \fill [brown, opacity = 0.3] (4,0) rectangle (6,4);
  \fill [brown, opacity = 0.3] (0,2) rectangle (2,6);
  \fill [brown, opacity = 0.3] (2,4) rectangle (6,6);
  \draw [line width=1pt] (0,0) rectangle (4,2);
  \draw [line width=1pt] (4,0) rectangle (6,4);
  \draw [line width=1pt] (0,2) rectangle (2,6);
  \draw [line width=1pt] (2,4) rectangle (6,6);
  \draw [line width=1pt, dotted] (4,0) -- (6,4);
  \draw [line width=1pt, dotted] (6,4) -- (2,6);
  \draw [line width=1pt, dotted] (2,6) -- (0,2);
  \draw [line width=1pt, dotted] (0,2) -- (4,0);
  \draw (1.8,-0.8) node[anchor=south west, scale=1.2] {$a$};
  \draw (-0.8,1) node[anchor=south west, scale=1.2] {$b$};
  \draw (1.8,1) node[anchor=south west, scale=1.2] {$c$};  
  \draw [decorate, decoration = {brace}]
      (2.05,2.1) --  (3.95,2.1);
  \draw (3,2.2) node[anchor=south] {$a-b$};
\end{tikzpicture}
\end{center}
\caption{Four Bricks and Pythagoras}\label{Abb4B}
\end{figure}
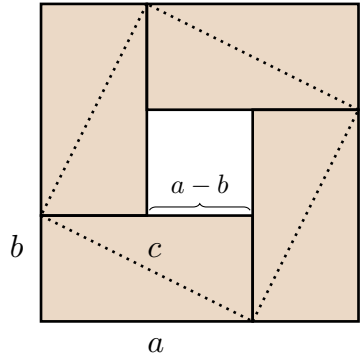

Drawing the diagonals of the four rectangles in the Four Bricks
diagram yields a square of side length $c$ (see Fig.~\ref{Abb4B}). The
angles at the vertices of this inner square must be right angles,
since each rectangle is rotated by $90^\circ$ with respect to its
neighbor.

The area of this square is clearly
$$ c^2 = (a - b)^2 + 2ab = a^2 + b^2, $$
which immediately yields the

\begin{thm}[Pythagorean Theorem]
  In a right triangle with legs $a$ and $b$ and hypotenuse $c$, we have
  $$ a^2 + b^2 = c^2. $$
\end{thm}

\subsection{Pythagorean Triples}

If we compute the area of the four-bricks figure in two different ways,
then we obtain the equation
$$ (a - b)^2 + 4ab = (a + b)^2. $$
Since $(a - b)^2 \ge 0$, it follows, as we have already seen, that
$(a + b)^2 \ge 4ab$, which is precisely the inequality in Theorem~\ref{SAGM}.
Now substitute $a = m^2$ and $b = n^2$ into this equation. This yields
\begin{equation}\label{PT}
  (m^2 - n^2)^2 + (2mn)^2 = (m^2 + n^2)^2.
\end{equation}
Thus we have the following result:

\begin{prop}
  For all natural numbers $m$ and $n$, equation (\ref{PT}) holds.
  In particular,
  $$ (m^2 - n^2,\ 2mn,\ m^2 + n^2) $$
  is a Pythagorean triple.
\end{prop}

\subsection{Quadratic Equations}

The solution of the quadratic equation
\begin{equation}\label{Eqpq}
  X^2 - pX + q = 0,
\end{equation}
or, equivalently by Vieta's formulas, of the system
$$ x + y = p, \quad xy = q, $$
where $x$ and $y$ are the solutions of (\ref{Eqpq}),
can also be visualized in an elegant manner with the four-bricks figure.

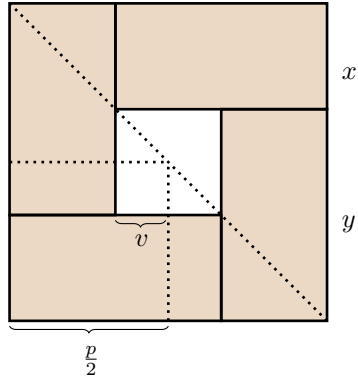
\begin{figure}[ht!]
\begin{center}
\begin{tikzpicture}[scale=0.7]
  \fill [brown, opacity = 0.3] (0,0) rectangle (4,2);
  \fill [brown, opacity = 0.3] (4,0) rectangle (6,4);
  \fill [brown, opacity = 0.3] (0,2) rectangle (2,6);
  \fill [brown, opacity = 0.3] (2,4) rectangle (6,6);
  \draw [line width=1pt] (0,0) rectangle (4,2);
  \draw [line width=1pt] (4,0) rectangle (6,4);
  \draw [line width=1pt] (0,2) rectangle (2,6);
  \draw [line width=1pt] (2,4) rectangle (6,6);
  \draw [line width=1pt, dotted] (0,6) -- (6,0);
  \draw [line width=1pt, dotted] (0,3) -- (3,3);
  \draw [line width=1pt, dotted] (3,0) -- (3,3);
  \draw [decorate, decoration = {brace}]
      (3,-0.1) --  (0,-0.1);
  \draw (6.1,1.5) node[anchor=south west] {$y$};
  \draw (6.1,4.4) node[anchor=south west] {$x$};
  \draw (1.2,-1.2) node[anchor=south west] {$\frac{p}{2}$};
  \draw [decorate, decoration = {brace}]
      (2.95,1.9) --  (2,1.9);
  \draw (2.5,1.8) node[anchor=north] {$v$};
\end{tikzpicture}
\end{center}
\caption{Four Bricks and quadratic equations}
\end{figure}

In the figure above it is apparent that
$$ \frac{p^2}{4} = xy + v^2, $$
and since $xy = q$, it follows that
$$ v^2 = \frac{p^2}{4} - q. $$
On the other hand we have
$$ x = \frac{p}{2} + v \quad \text{and} \quad y = \frac{p}{2} - v. $$
With this we obtain the

\begin{prop}
  The solutions $x$ and $y$ of the system of equations
  $$   x + y = p, \quad xy = q   $$
  are given by
  $$  x = \frac{p}{2} - \sqrt{\frac{p^2}{4} - q}, \qquad
      y = \frac{p}{2} + \sqrt{\frac{p^2}{4} - q}.  $$
\end{prop}

\section{Counting Bricks}

As the Babylonians must have known a connection between bisecting
trapezoids with rational sides and Pythagorean triples, and since
Euclid's proof given in Section~\ref{S1} does little to explain this
connection, we will look for other ways of deriving the Babylonian equation.

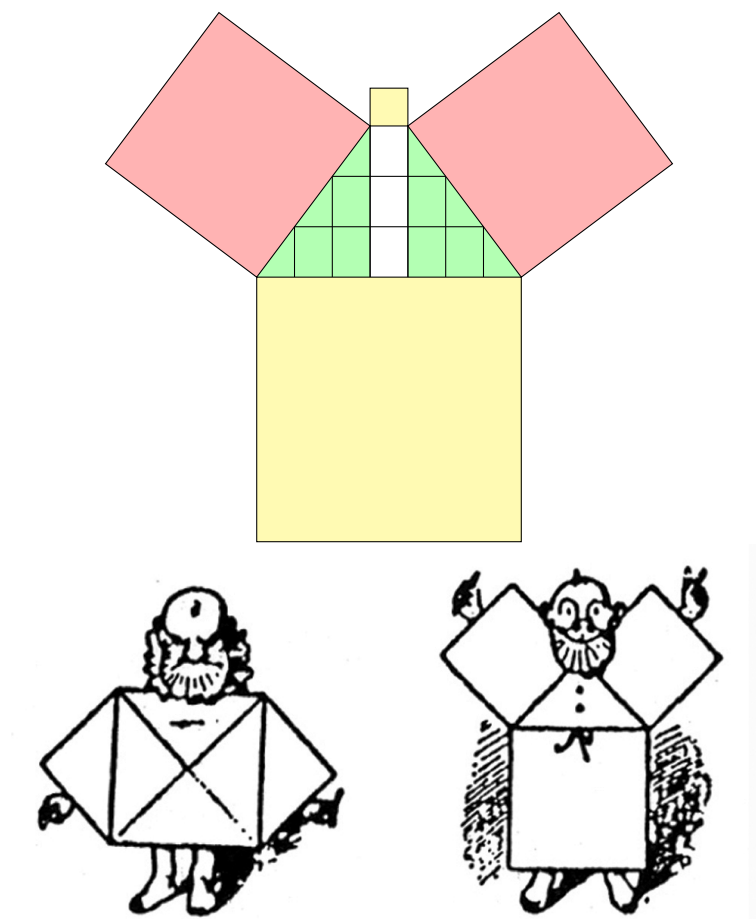
\begin{figure}[p!]
\begin{center}
\begin{tikzpicture}[scale=0.5]
  \fill [green, opacity = 0.3] (0,0) -- (3,0) -- (3,4) -- cycle;
  \fill [green, opacity = 0.3] (4,0) -- (7,0) -- (4,4) -- cycle;  
  \fill [red, opacity=0.3] (-4,3) -- (0,0) -- (3,4) -- (-1,7) -- cycle;
  \fill [red, opacity=0.3] (7,0) -- (11,3) -- (8,7) -- (4,4) -- cycle;
  \fill [yellow, opacity=0.3] (3,4) rectangle (4,5);
  \fill [yellow, opacity=0.3] (0,0) rectangle (7,-7);
  \draw (-4,3) -- (0,0) -- (3,4) -- (-1,7) -- cycle;
  \draw (7,0) -- (11,3) -- (8,7) -- (4,4) -- cycle;
  \draw (3,4) rectangle (4,5);
  \draw (0,0) rectangle (7,-7);
  \draw (3,0) -- (3,4);
  \draw (4,0) -- (4,4);
  \begin{scope}
    \clip (0,0) -- (7,0) -- (4,4) -- (3,4) -- cycle;
    \foreach \x in {1,2,3,4,5,6} \draw (\x,0) -- (\x,4);
    \foreach \x in {1,2} \draw (0,4*\x/3) -- (7,4*\x/3);    
  \end{scope}  
\end{tikzpicture}

\includegraphics[width=10cm]{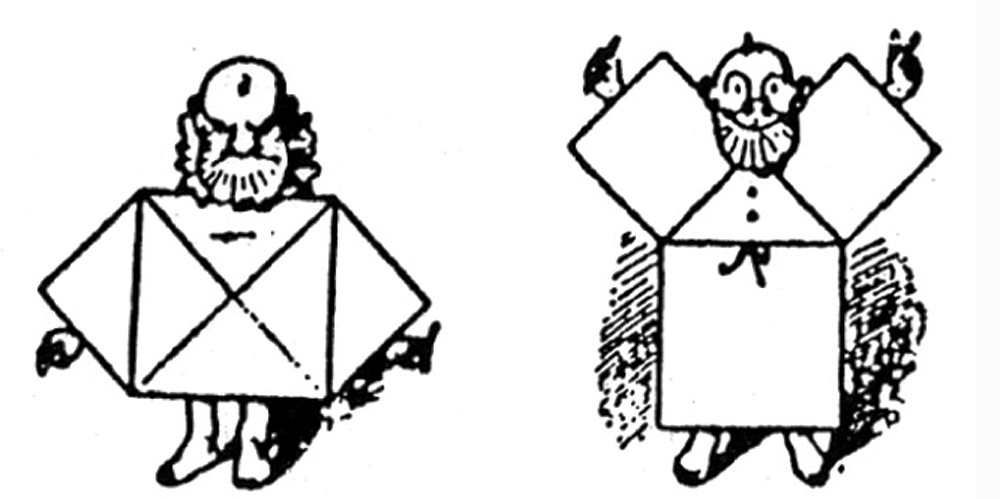}

\end{center}
\caption{Top: The ``Babylonian Pythagoras'' based on the
  Babylonian triple $(1,5,7)$ inspired by the cartoon below. Bottom:
  Proof without words of the Pythagorean theorem for equilateral right
  angled triangles in ``Fliegende Bl \"atter'' (vol. 49 (1868), p. 52),
  titled ``Pythagoras before and after the discovery of his theorem''.}
\label{FigBT}
\end{figure}
\clearpage{}

We will now explain how to generate the Babylonian trapezoid $(1,5,7)$
from the Pythagorean triangle with sides $(3,4,5)$.  The two green
Pythagorean triangles in Fig.~\ref{FigBT} give us the relation
$3^2 + 4^2 = 5^2$; the two red squares thus have the area
$2 \cdot 5^2 = 2 \cdot 3^2 + 2 \cdot 4^2$. The identity (\ref{E4B4})
already used above shows that 
$$ 5^2 + 5^2 = (4-3)^2 + (4+3)^2 = 1^2 + 7^2, $$
hence the red and the yellow squares have the same area.

If we rearrange the trapezoid into a square consisting of a green
rectangle plus a white stripe, we see that its area is equal to $4^2$.
If we tile the trapezoid as indicated, we see that the lower
transversal halves the area of the trapezoid. The length of this
transversal is obviously $1 + \frac{2}{3}(7-1) = 5$.

\section{Division of Trapezoids using Four Bricks}

As a first geometric interpretation of this construction we
consider the figure with the ``four bricks'' that we have already
used for justifying the Pythagorean rule\footnote{The basic idea
  of this proof can be found in \cite{BBS}.}:

\begin{figure}[ht!]
\begin{center}
\begin{tikzpicture}[scale=0.6]
  \draw (0,0) rectangle (7,7);
  \draw (0,4) -- (3,0) -- (7,3) -- (4,7) -- cycle;
  \draw (0,4) -- (4,4);
  \draw (3,4) -- (3,0);
  \draw (3,3) -- (7,3);
  \draw (4,3) -- (4,7);
  \draw (1.5,-0.7) node[anchor=south] {$a$};
  \draw (  5,-0.7) node[anchor=south] {$b$};
\end{tikzpicture} \quad 
\begin{tikzpicture}[scale=0.6]
  \fill [brown, opacity=0.2] (1,1) rectangle (6,6);
  \draw (1,1) rectangle (6,6);
  \draw (0,0) rectangle (7,7);
  \draw [red] (0,0) -- (3,3) -- (4,3) -- (7,0) -- cycle;
  \draw[dashed] (0,4) -- (3,0) -- (7,3) -- (4,7) -- cycle;
  \draw[dashed] (0,4) -- (4,4);
  \draw[dashed] (3,4) -- (3,0);
  \draw[dashed] (3,3) -- (7,3);
  \draw[dashed] (4,3) -- (4,7);
  \draw (1.5,-0.7) node[anchor=south] {$a$};
  \draw (  5,-0.7) node[anchor=south] {$b$};
\end{tikzpicture} 
\end{center}
\caption{Four Bricks and identities}
\end{figure}
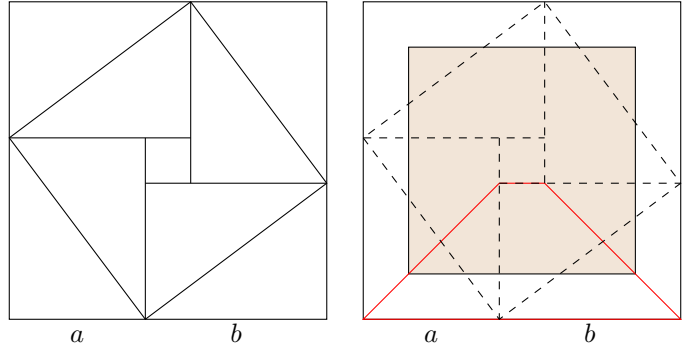

Now we turn the square in the middle until its sides are parallel to the
sides of the outer square (see the figure on the right). The trapezoid
appearing in this figure has base angles of $45^\circ$, hence the corners
of the square lie on the legs of the trapezoid, and the argument above
shows that the transversal dividing the trapezoid into two trapezoids
of the same area has length $c$.

\begin{figure}[ht!]
\begin{center}
\begin{tikzpicture}[scale=0.6]
  \fill [green, opacity=0.2] (0,0) -- (3,0) -- (3,3);
  \fill [red,   opacity=0.2] (3,0) -- (7,0) -- (4,3) -- (3,3) -- cycle;
  \draw (1,1) rectangle (6,6);
  \draw (0,0) rectangle (7,7);
  \draw [red] (0,0) -- (3,3) -- (4,3) -- (7,0) -- cycle;
  \draw[dashed] (0,4) -- (3,0) -- (7,3) -- (4,7) -- cycle;
  \draw[dashed] (0,4) -- (4,4);
  \draw[dashed] (3,4) -- (3,0);
  \draw[dashed] (3,3) -- (7,3);
  \draw[dashed] (4,3) -- (4,7);
  \draw (1.5,-0.75) node[anchor=south] {$a$};
  \draw (  5,-0.75) node[anchor=south] {$b$};
\end{tikzpicture} \quad 
\begin{tikzpicture}[scale=0.6]
  \fill [red,   opacity=0.2] (3,0) -- (7,0) -- (4,3) -- (3,3) -- cycle;
  \fill [green, opacity=0.2] (7,0) -- (7,3) -- (4,3) -- cycle;
  \draw (1,1) rectangle (6,6);
  \draw (0,0) rectangle (7,7);
  \draw [red] (0,0) -- (3,3) -- (4,3) -- (7,0) -- cycle;
  \draw[dashed] (0,4) -- (3,0) -- (7,3) -- (4,7) -- cycle;
  \draw[dashed] (0,4) -- (4,4);
  \draw[dashed] (3,4) -- (3,0);
  \draw[dashed] (3,3) -- (7,3);
  \draw[dashed] (4,3) -- (4,7);
  \draw (1.5,-0.75) node[anchor=south] {$a$};
  \draw (  5,-0.75) node[anchor=south] {$b$};
  \draw [decorate, decoration={brace}] (6.05,7.1) -- (7,7.1);
  \draw (6.5,7.2) node[anchor=south] {$\frac{a+b-c}2$};
  \draw [decorate, decoration={brace}] (1.,1.1) -- (6,1.1);
  \draw (3.5,1.3) node[anchor=south] {$c$};
\end{tikzpicture} 
\end{center}
\caption{Bisecting Trapezoids}
\end{figure}
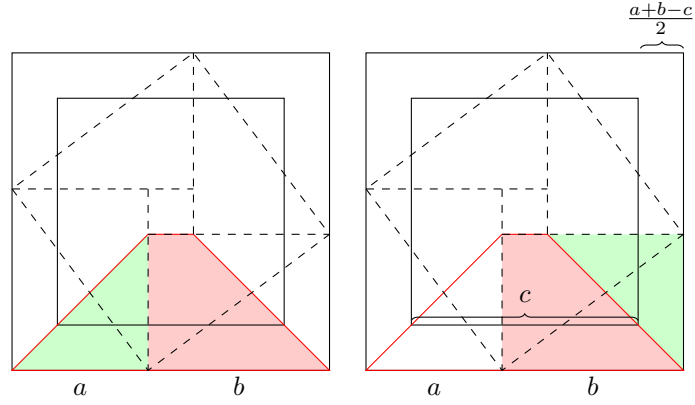
  
Given a right angled triangle with sides $(a, b, c)$ with $a < b$, 
the red trapezoid has parallel sides $a' = b + a$ and $c' = b - a$,
and its bisecting transversal has length $c$, giving rise to the
Babylonian triple $(b+a,c,b-a)$.

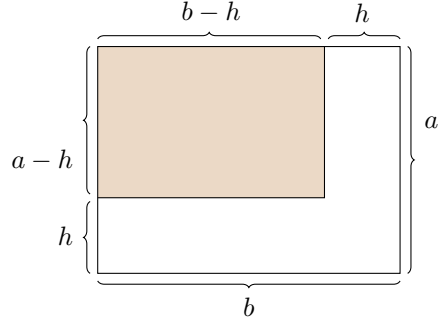
\begin{figure}[ht!]
\begin{center}
  \begin{tikzpicture}
    \fill [brown, opacity=0.3] (0,1) rectangle (3,3);
    \draw (0,0) rectangle (4,3);
    \draw (0,1) -- (3,1) -- (3,3);
    \draw [decorate, decoration={brace}] (-0.1,1.05) -- (-0.1,3);
    \draw (-0.2,1.5) node[anchor=east] {$a-h$};
    \draw [decorate, decoration={brace}] (-0.1,0) -- (-0.1,0.95);
    \draw (-0.2,0.5) node[anchor=east] {$h$};
    \draw [decorate, decoration={brace}] (0,3.1) -- (2.95,3.1);
    \draw (1.5,3.2) node[anchor=south] {$b-h$};
    \draw [decorate, decoration={brace}] (3.05,3.1) -- (4,3.1);
    \draw (3.5,3.2) node[anchor=south] {$h$};
    \draw [decorate, decoration={brace}] (4.1,3) -- (4.1,0);
    \draw (4.2,2) node[anchor=west] {$a$};
    \draw [decorate, decoration={brace}] (4,-0.1) -- (0,-0.1);
    \draw (2,-0.2) node[anchor=north] {$b$};
  \end{tikzpicture} 
  \caption{Dividing a rectangle}\label{AbbHR}
\end{center}
\end{figure}

A second but more involved argument is the following.  If we move the
green Pythagorean triangle to the right, we obtain a rectangle with
area $ab$ (see Fig.~\ref{AbbHR}). The brown rectangle in the figure
below has area $A = ab$. The smaller rectangle in the left upper
corner has half its area; if the height of the transversal is $h$,
then we must have
$$ (a-h)(b-h) = \frac{ab}2, \quad \text{i.e.,} \quad
   h^2 - (a+b)h + \frac{ab}2 = 0. $$
Solving this quadratic equation using $a^2 + b^2 = c^2$ we find that 
$$ h = \frac{a+b-c}2. $$

\begin{figure}[ht!]
  \begin{center}   
  \begin{tikzpicture}[scale=0.7]
  \fill[green, opacity=0.3] (0,0) -- (3,0) -- (3,3) -- cycle;
  \fill[green, opacity=0.3] (7,0) -- (4,0) -- (4,3) -- cycle;
  \fill[red,   opacity=0.3] (3,0) -- (4,0) -- (4,3) -- (3,3) -- cycle;
  \draw (0,0) -- (3,3) -- (4,3) -- (7,0) -- cycle;
  \draw[dotted] (3,0) -- (3,3);
  \draw[dotted] (4,0) -- (4,3);
  \draw[dashed] (1,1) -- (6,1);
  \end{tikzpicture} \qquad     
  \begin{tikzpicture}[scale=0.7]
  \fill[green, opacity=0.3] (0,0) -- (3,0) -- (3,4) -- cycle;
  \fill[green, opacity=0.3] (7,0) -- (4,0) -- (4,4) -- cycle;
  \fill[red,   opacity=0.3] (3,0) -- (4,0) -- (4,4) -- (3,4) -- cycle;
  \draw (0,0) -- (3,4) -- (4,4) -- (7,0) -- cycle;
  \draw[dotted] (3,0) -- (3,4);
  \draw[dotted] (4,0) -- (4,4);
  \draw[dashed] (1,{4/3}) -- (6,{4/3});
  \end{tikzpicture}
  \caption{Trapezoid constructed from Pythagorean triangles.}\label{AbbT2}
\end{center}
\end{figure}
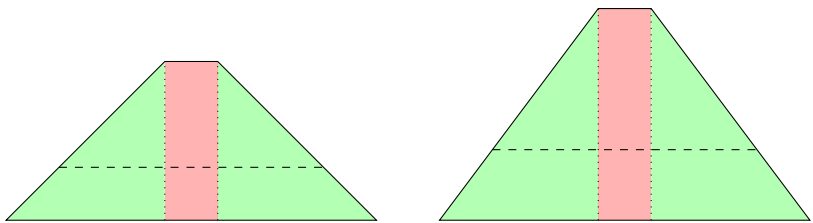

The Four-Bricks diagram built from the Pythagorean triple $(a, b, c)$
with $a < b$ contains two distinct trapezoids with parallel sides
$a' = b + a$ and $c' = b-a$ with bisecting transversal of length $c$:
the one showing up in the first proof above, which has height $a$
and a base angle of $45^\circ$ (Fig.~\ref{AbbT2} left),
and the one built from two copies of the triangle with sides $(a, b, c)$
and the rectangle with length $b-a$ and height $b$ (Fig.~\ref{AbbT2} right).
The right diagram displays the right angled triangles $(a,b,c)$
from which it is constructed, and this diagram shows up in Fig.~\ref{FigBT}.
In the left diagram, the right-angled triangle is hidden, but the
derivation of the Babylonian triple $(a', b', c')$ is quite simple.

These constructions provide us with a Babylonian trapezoid for every
Pythagorean triangle (see Fig.~\ref{AbbT2} right): Insert a rectangle
with base $b-a$ and height $b$ between the two Pythagorean triangles
$(a, b, c)$; then the trapezoid with parallel sides $a+b$ and $b-a$
has a transversal of length $c$ dividing the trapezoid into trapezoids
with equal areas.


\begin{thebibliography}{99}

\bibitem{Bonet} A. Bonet,
  \url{https://www.researchgate.net/profile/Aldo-Bonet}
  %
  
\bibitem{BBS} L. Brack-Bernsen, O. Schmidt,
{\em Bisectable trapezia in Babylonian mathematics},
Centaurus {\bf 33} (1990), 1--38
%

\bibitem{Caveing6} M. Caveing,
  {\em La tablette babylonienne AO 17264 du mus\'ee du louvre et
    le probl\`eme des six fr\`eres},
Hist. Math. {\bf 12} (1985), 6--24
%

\bibitem{Cid} J. A. Cid,
  {\em Una Figura, Dos Identidades, \textexclamdown{}Seis Teoremas!},
  Gaceta de la RSME {\bf 20} (2017), 398
  %
  
\bibitem{DamPyth} P. Damerow,
{\em Kannten die Babylonier den Satz des Pythagoras?},
Antrittsvorlesung 1994; in \cite[219--310]{HoyDam}
%

\bibitem{ED} R. C. Archibald, 
  {\em Euclid's book on divisions of figures},
  Cambridge University Press 1915
  %
  
\bibitem{Frib} J. Friberg,
{\em A remarkable collection of Babylonian mathematical texts},
Springer-Verlag 2007
%

\bibitem{FribT} J. Friberg,
{\em Amazing Traces of a Babylonian Origin in Greek Mathematics},
World Scientific 2007
%

\bibitem{Gandz} S. Gandz,
{\em Studies in Babylonian Mathematics I: Indeterminate Analysis 
      in Babylonian Mathematics},
Osiris {\bf 8} (1948), 12--40
%

\bibitem{GoncTH} C. Gon\c{c}alves,
{\em Mathematical Tablets from Tell Harmal},
Springer-Verlag 2015
%

\bibitem{Hoy} J. H\o{}yrup,
  {\em Algebra in Cuneiform. Introduction to an Old Babylonian
    geometrical technique}, MPI Hist. Sci. 2017
  %
  
\bibitem{HoyDam} J. H\o{}yrup, P. Damerow (ed.),
{\em Changing Views on Ancient Near Eastern Mathematics},
D. Reimer Verlag Berlin, 2001
%

\bibitem{LemMI} F. Lemmermeyer,
{\em Zur Zahlengeometrie der Babylonier},
Mathematik-Information 68 (2018), 15--30
%

\bibitem{LemmBA} F. Lemmermeyer,
  {\em Mathematik \`a la carte. Babylonische Algebra},
  Springer Spektrum 2022
  %
  
  
\bibitem{Lem4000} F. Lemmermeyer,
  {\em 4000 Jahre Zahlentheorie. Geschichte - Kulturen - Menschen I.
    Von Babel bis Abel}, Springer Spektrum  2023; English Translation
  Birkh\"auser 2026, to appear
  %

\bibitem{Neugebauer} O. Neugebauer,  
{\em Mathematische Keilschrifttexte} I, Springer 1935, 
II 1935, III 1937
%

 \bibitem{Polk} M. Polk, A. M. H. Schuster (ed.),
  {\em The looting of the Iraq Museum, Baghdad. The lost legacy of
    Ancient Mesopotamia}, New York 2005
  %
  
  
  
\end{thebibliography}
\end{document}